# FREQUENCY OF A DIGIT IN THE REPRESENTATION OF A NUMBER AND THE ASYMPTOTIC MEAN VALUE OF THE DIGITS

**S. O. Klymchuk,**[1] **O. P. Makarchuk,**[2] **and M. V. Prats'ovytyi**[1,2]   UDC 511.72+517+519.21

We study the relationship between the frequency of a ternary digit in a number and the asymptotic mean value of the digits. The conditions for the existence of the asymptotic mean of digits in a ternary number are established. We indicate an infinite everywhere dense set of numbers without frequency of digits but with the asymptotic mean of the digits.

## 1. Introduction

It is known that, for *any real number* $x \in [0;1]$, *there exists a sequence* $(\alpha_n)$ *such that*

$$\alpha_n \in \mathcal{A} = \{0, 1, \ldots, s-1\}$$

*and*

$$x = \frac{\alpha_1}{s} + \frac{\alpha_2}{s^2} + \ldots + \frac{\alpha_n}{s^n} + \ldots \equiv \Delta^s_{\alpha_1\alpha_2\ldots\alpha_n\ldots}.$$

The last relation is called the $s$-*nary representation* and $\alpha_k = \alpha_k(x)$ is called the $k$ th $s$-*nary digit of the number* $x$. However, the $k$ th digit of the number is, generally speaking, ill-defined as its function because the following equality is true:

$$\Delta^s_{c_1\ldots c_{k-1}c_k(0)} = \Delta^s_{c_1\ldots c_{k-1}[c_k-1](s-1)},$$

where $(i)$ is the period in the representation of the number. The numbers of this kind are called $s$-*nary-rational*; they have exactly two $s$-nary representations. The other numbers have only one representation; they are called $s$-*nary-irrational*. For the correct definition of the $k$ th digit, it suffices to agree to use only the first $s$-nary representation, namely, the representation with period $(0)$. As a result, the function $\alpha_k(x)$ becomes well defined on $[0;1]$.

Let $N_i(x,k)$ be the number of digits $i \in \mathcal{A} = \{0, \ldots, s-1\}$ in the $s$-nary representation $\Delta^s_{\alpha_1\alpha_2\ldots\alpha_k\ldots}$ of the number $x \in [0;1]$ up to the $k$ th place inclusively, i.e.,

$$N_i(x,k) = \#\{j \colon \alpha_j(x) = i, \ j \leqslant k\}.$$

The number $v_i^{(n)} \equiv k^{-1}N_i(x,k)$ is the relative frequency of the digit $i$ in the $s$-nary representation of the real number $x$.

**Definition 1.** *The frequency (asymptotic frequency) of the digit* $i$ *in the* $s$-*nary representation of a number* $x \in [0;1]$ *is defined as the limit*

---

[1] Institute of Mathematics, Ukrainian National Academy of Sciences, Kyiv, Ukraine.
[2] Drahomanov National Pedagogic University, Kyiv, Ukraine.

$$\nu_i(x) = \lim_{k \to \infty} \frac{N_i(x, k)}{k}$$

*if it exists.*

It is clear that the function of frequency $\nu_i(x)$ of the digit $i$ in the $s$-nary representation of a number $x \in [0; 1]$ is well defined for $s$-nary-irrational numbers; for the $s$-nary-rational numbers, this function becomes well defined only after the above-mentioned agreement concerning the use solely of the representation with period (0).

It is known that the function of frequency of an $s$-nary digit of the number takes all values from $[0, 1]$; at rational points from $[0, 1]$, it is defined and takes rational values; it is not defined on an everywhere dense continual set of numbers from $[0, 1]$ of Lebesgue measure zero. This notion is useful for the investigation of pure distributions of random variables whose $s$-nary digits are random [7, 8].

The number $r_n \equiv \frac{1}{n} \sum_{i=1}^{n} \alpha_i(x)$ is called the relative mean of digits for the number $x$. Since

$$r_n(x) = \frac{N_1(x, k)}{k} + \frac{2N_2(x, k)}{k}, \qquad \text{we have} \qquad 0 \leqslant r_n(x) \leqslant 2.$$

**Definition 2.** *If there exist a limit*

$$\lim_{n \to \infty} \frac{1}{n} \sum_{i=1}^{n} \alpha_i(x) = r(x),$$

*where $\alpha_i$ are the $s$-nary digits of the number $x \in [0; 1]$, then its value (the number $r(x)$) is called the asymptotic mean (or simply mean) value for the digits of the number $x$.*

The asymptotic mean value of digits is a certain analog of the frequency of an $s$-nary digit of the number. Moreover, for $s = 2$, the equality $r(x) = \nu_1(x)$ is true.

We are interested in topological metric properties of the sets of numbers with *a priori* specified *asymptotic mean value of digits*, i.e., sets of the form

$$S_\theta = \left\{ x : \lim_{n \to \infty} \frac{1}{n} \sum_{i=1}^{n} \alpha_i(x) = \theta \geqslant 0 \right\},$$

where the constant $\theta$ is an *a priori* given parameter from the segment $[0, s-1]$.

Note that, for $\theta > s - 1$, the set $S_\theta$ is empty.

It is easy to see that the set $S_\theta$ is, to a certain extent, similar to the Besicovitch–Eggleston set

$$E\big[\tau_0, \tau_1, \ldots, \tau_{s-1}\big] = \big\{ x : \nu_i(x) = \tau_i, \ i = \overline{0, s-1} \big\}.$$

In what follows, we establish the relationship between these sets.

## 2. Relationship between the Frequency and the Mean of Digits

**Lemma 1.** *If the $s$-nary representation of a number $x$ has frequencies of all digits $\nu_0, \nu_1, \ldots, \nu_{s-1}$, then it has the asymptotic value of digits $r(x)$ and, in addition,*

$$r(x) = \nu_1(x) + 2\nu_2(x) + \ldots + (s-1)\nu_{s-1}(x).$$

*Proof.* Since

$$\frac{1}{n}\sum_{i=1}^{n}\alpha_i(x) = \frac{0\cdot N_0(x,n)}{n} + \frac{1\cdot N_1(x,n)}{n} + \ldots + \frac{(s-1)\cdot N_{s-1}(x,n)}{n}$$

$$= \frac{N_1(x,n)}{n} + 2\cdot\frac{N_2(x,n)}{n} + \ldots + (s-1)\cdot\frac{N_{s-1}(x,n)}{n},$$

we have

$$\lim_{n\to\infty}\frac{1}{n}\sum_{i=1}^{n}\alpha_i(x) = \lim_{n\to\infty}\left(\frac{N_1(x,n)}{n} + 2\cdot\frac{N_2(x,n)}{n} + \ldots + (s-1)\cdot\frac{N_{s-1}(x,n)}{n}\right)$$

$$= \lim_{n\to\infty}\frac{N_1(x,n)}{n} + 2\lim_{n\to\infty}\frac{N_2(x,n)}{n} + \ldots + (s-1)\lim_{n\to\infty}\frac{N_{s-1}(x,n)}{n}$$

if the last limits exist. By using the last equality, we obtain

$$r(x) = \nu_1(x) + 2\nu_2(x) + \ldots + (s-1)\nu_{s-1}(x).$$

*Corollary 1.* If $\theta = \tau_1 + 2\tau_2 + \ldots + (s-1)\tau_{s-1}$, then the following inclusion is true:

$$E[\tau_0, \tau_1, \ldots, \tau_{s-1}] \subset S_\theta.$$

**Theorem 1.** *The property of the number $r(x) = \dfrac{s-1}{2}$ is normal, i.e., the set of real numbers from the segment $[0;1]$ without this property is a set of Lebesgue measure zero.*

*Proof.* According to the well-known Borel theorem [4], for almost all numbers from the segment $[0;1]$ (in a sense of Lebesgue measure), the frequencies of all digits in their $s$-nary representation exist and are equal to $s^{-1}$, i.e., the Besicovitch–Eggleston set $E\left[s^{-1}, s^{-1}, \ldots, s^{-1}\right]$ is a set of complete Lebesgue measure. Since

$$\theta = \frac{1}{s} + \frac{2}{s} + \frac{3}{s} + \ldots + \frac{s-1}{s} = \frac{s-1}{2},$$

we have $E\left[s^{-1}; s^{-1}; \ldots; s^{-1}\right] \subset S_\theta$. Therefore, $\lambda(S_\theta) = 1$.

**Theorem 2.** *If, for $s = 3$, there exist an asymptotic mean value of the digits $r(x)$ and at least one of the frequencies $\nu_0(x)$, $\nu_1(x)$, and $\nu_2(x)$, then the other two frequencies of digits for the number $x$ also exist. If at least one of the frequencies $\nu_j(x)$, $j \in \{0,1,2\}$, does not exist but $r(x)$ exists, then the other two frequencies of digits of the number $x$ also do not exist.*

*Proof.* Let $v_j^{(n)} = n^{-1}N_j(x,n)$ be the relative frequency of a ternary digit $j$ of a number $x$ and let $r_n(x) = \dfrac{1}{n}\sum_{j=1}^{n}\alpha_j(x)$ be the relative mean value of digits for the number $x$. Then we get the following system:

$$v_0^{(n)} + v_1^{(n)} + v_2^{(n)} = 1,$$

$$v_1^{(n)} + 2v_2^{(n)} = r_n.$$

If the limits $\lim_{n\to\infty} r_n$ and $\lim_{n\to\infty} v_j^{(n)}$, $j \in \{1, 2\}$, exist, then the second equality of the system implies that $\lim_{n\to\infty} v_i^{(n)}$, $i \in \{1, 2\}\setminus\{j\}$ and the first equality yields the limit $\lim_{n\to\infty} v_0^{(n)}$. Since the system can be rewritten in the form

$$v_2^{(n)} = r_n - 1 + v_0^{(n)},$$

$$v_1^{(n)} = 2 - 2v_0^{(n)} - r_n,$$

the existence of $\lim_{n\to\infty} v_0^{(n)}$ and $\lim_{n\to\infty} r_n$ implies the existence of the frequencies $\nu_j(x) = \lim_{n\to\infty} v_j^{(n)}$, $j = 1, 2$.

If the limit $\lim_{n\to\infty} r_n$ exists but the limit $\lim_{n\to\infty} v_0^{(n)}$ does not exist, then the limits $\lim_{n\to\infty} v_1^{(n)}$ and $\lim_{n\to\infty} v_2^{(n)}$ do not exist and, hence, the frequencies $\nu_1(x)$ and $\nu_2(x)$ also do not exist.

***Corollary 2.*** *If the ternary representation of a number has the asymptotic mean value $r(x)$, then the frequencies $\nu_i(x)$ $i \in \{0, 1, 2\}$, either simultaneously exist or simultaneously do not exist.*

## 3. Numbers with Preassigned Asymptotic Mean Value of the Digits

We now present an algorithm for the construction of a number with preassigned frequencies of the ternary digits.

Let $a$ and $b$ be positive real numbers such that $a + b \leqslant 1$ and, in addition, $\nu_0(x) = a$ and $\nu_1(x) = b$. By $(c_n)$, where $c_n = [n \cdot a]$, $n \in \mathbb{N}$, we denote an integer-valued sequence. Consider the difference

$$d_n = c_{n+1} - c_n = \big[(n+1) \cdot a\big] - [n \cdot a] = \big[[n \cdot a] + \{n \cdot a\} + a\big] - [n \cdot a] = \big[\{n \cdot a\} + a\big] \in \{0, 1\}.$$

Then $x_* = \Delta^3_{\alpha_2 \alpha_3 \ldots \alpha_n \alpha_{n+1} \ldots}$, where

$$\alpha_{n+1} = \begin{cases} \beta_n & \text{if } d_n = 0, \\ 0 & \text{if } d_n = 1. \end{cases}$$

We introduce an integer-valued sequence $(c'_n)$, $c'_n = [n \cdot b]$, $n \in \mathbb{N}$, and consider the difference $d'_n = c'_{n+1} - c'_n \in \{0, 1\}$. As a result, we obtain

$$\beta_n = \begin{cases} 1 & \text{for } d'_n = 0, \\ 2 & \text{for } d'_n = 1. \end{cases}$$

Since

$$a - \frac{1}{n} = \frac{n \cdot a - 1}{n} < \frac{N_0(x_*, n)}{n} = \frac{[n \cdot a]}{n} \leqslant \frac{n \cdot a}{n} = a,$$

$$b - \frac{1}{n} = \frac{n \cdot b - 1}{n} < \frac{N_1(x_*, n)}{n} = \frac{[n \cdot b]}{n} \leqslant \frac{n \cdot b}{n} = b,$$

we find

$$\nu_0(x_*) = \lim_{n \to \infty} \frac{N_0(x_*, n)}{n} = a, \qquad \nu_1(x_*) = \lim_{n \to \infty} \frac{N_1(x_*, n)}{n} = b.$$

## 4. Numbers for Which with the Asymptotic Mean Value of Digits Exists but the Frequencies of Digits Do Not Exist

We now construct a number for which the asymptotic frequencies of digits do not exist but the asymptotic mean value of digits exists and is equal to a certain given number $\theta$, i.e., $r(x) = \theta$.

First, we show that this cannot be done for $\theta = 0$ and $\theta = 2$.

Let $\theta = 0$, i.e., $\lim_{n \to \infty} r_n = 0$. Since $v_1^{(n)} + 2v_2^{(n)} \geqslant v_i^{(n)} \geqslant 0$ for $i \in \{1, 2\}$, we have $\lim_{n \to \infty} v_i^{(n)} = 0$. Now let $\theta = 2$. Since

$$r_n = v_1^{(n)} + 2v_2^{(n)} = v_1^{(n)} + v_2^{(n)} + v_2^{(n)} = 1 - v_0^{(n)} + v_2^{(n)},$$

we conclude that

$$1 \geqslant v_2^{(n)} = r_n - 1 + v_0^{(n)} \geqslant r_n - 1.$$

Further, since $\lim_{n \to \infty} r_n = 2$, we get $\lim_{n \to \infty} v_2^{(n)} = 1$. Hence,

$$0 \leqslant v_i^{(n)} = 1 - v_2^{(n)} - v_{1-i}^{(n)} \leqslant 1 - v_2^{(n)},$$

where $i \in \{0, 1\}$. However, $\lim_{n \to \infty} 1 - v_2^{(n)} = 0$ and, therefore, $\lim_{n \to \infty} v_i^{(n)} = 0$, $i \in \{0, 1\}$.

Let $\theta$ be a given number from $(0; 2)$. We now construct a representation of the number $x^*$ with the asymptotic mean value of digits $\theta$ but without any frequency of ternary digits. To this end, we consider the following representation of $x^*$:

$$x^* = \Delta^3_{\underbrace{0\ldots0}_{a_{11}}\underbrace{1\ldots1}_{a_{12}}\underbrace{2\ldots2}_{a_{13}}\underbrace{0\ldots0}_{a_{21}}\underbrace{1\ldots1}_{a_{22}}\underbrace{2\ldots2}_{a_{23}}\ldots\underbrace{0\ldots0}_{a_{k1}}\underbrace{1\ldots1}_{a_{k2}}\underbrace{2\ldots2}_{a_{k3}}\ldots},$$

$$\underbrace{\phantom{a_{11}\ a_{12}\ a_{13}}}_{\text{1 st block}} \quad \underbrace{\phantom{a_{21}\ a_{22}\ a_{23}}}_{\text{2 nd block}} \quad \underbrace{\phantom{a_{k1}\ a_{k2}\ a_{k3}}}_{k\text{ th block}}$$

where $a_{ij} \in N$, and establish conditions for the matrix $\|a_{ij}\|$ under which the limit

$$\lim_{k \to \infty} \frac{a_{11} + a_{21} + \ldots + a_{k1}}{a_{11} + a_{12} + a_{13} + \ldots + a_{k1} + a_{k2} + a_{k3}}$$

does not exist and, hence, $\nu_0(x)$ also does not exist but

$$r(x) = \lim_{k \to \infty} \frac{0 \cdot (a_{11} + \ldots + a_{k1}) + 1 \cdot (a_{12} + \ldots + a_{k2}) + 2 \cdot (a_{13} + \ldots + a_{k3})}{a_{11} + a_{12} + a_{13} + \ldots + a_{k1} + a_{k2} + a_{k3}} = \theta$$

exists.

To this end, we prove two auxiliary statements.

**Lemma 2.** *For any natural numbers $k$ and $n$ $(k < n)$ and any real number $x$, the following equality holds:*

$$\lim_{n \to \infty} \frac{[kx] + [(k+1)x] + \ldots + [nx]}{\frac{n(n+1)}{2}} = x,$$

*where $[kx]$ is the integral part of the number $kx$.*

*Proof.* According to properties of the integral part of the number, we get $y - 1 < [y] \leq y$. This yields

(i) $$\frac{[kx] + [(k+1)x] + \ldots + [nx]}{\frac{n(n+1)}{2}} \leq \frac{kx + (k+1)x + \ldots + nx}{\frac{n(n+1)}{2}}$$

$$= \frac{x\bigl(1 + 2 + \ldots + n - (1 + 2 + \ldots + k - 1)\bigr)}{\frac{n(n+1)}{2}}$$

$$= \frac{x\frac{n(n+1)}{2} - x\frac{(k-1)k}{2}}{\frac{n(n+1)}{2}} \to x \quad \text{as} \quad n \to \infty;$$

(ii) $$\frac{[kx] + [(k+1)x] + \ldots + [nx]}{\frac{n(n+1)}{2}} > \frac{kx - 1 + (k+1)x - 1 + \ldots + nx - 1}{\frac{n(n+1)}{2}}$$

$$= \frac{kx + (k+1)x + \ldots + nx - (n - k + 1)}{\frac{n(n+1)}{2}}$$

$$= \frac{x\frac{n(n+1)}{2} - x\frac{(k-1)k}{2} - (n - k + 1)}{\frac{n(n+1)}{2}} \to x \quad \text{as} \quad n \to \infty.$$

Therefore,

$$\lim_{n \to \infty} \frac{[kx] + [(k+1)x] + \ldots + [nx]}{\frac{n(n+1)}{2}} = x.$$

We choose two arbitrary real numbers $x_1$ and $x_2$ and a real number $\varepsilon > 0$ such that $0 < x_1 < x_2$ and $\varepsilon < \frac{x_2 - x_1}{2}$. The last inequality is equivalent to the inequality $x_1 + \varepsilon < x_2 - \varepsilon$. We construct a sequence $(\alpha_n)$ (also denoted by $\tilde{\alpha}_n\{x_1, x_2\}$) each term of which takes one of two possible values $(x_1$ or $x_2)$ according to the rule:

Let $n_1$ be the least natural number for which the following inequality holds:

$$\frac{[1 \cdot x_1] + [2 \cdot x_1] + \ldots + [n_1 \cdot x_1]}{\frac{n_1(n_1 + 1)}{2}} < x_1 + \varepsilon$$

By Lemma 2, this number exists. Then we set $\alpha_1 = \alpha_2 = \ldots = \alpha_{n_1} = x_1$.

Let $n_2$ be the least natural number for which the inequality

$$\frac{[1 \cdot x_1] + \ldots + [n_1 \cdot x_1] + [(n_1 + 1) \cdot x_2] + \ldots + [n_2 \cdot x_2]}{\frac{n_2(n_2 + 1)}{2}} > x_2 - \varepsilon$$

is true. By Lemma 2, this number exists. Then we set $\alpha_{n_1+1} = \alpha_{n_1+2} = \ldots = \alpha_{n_2} = x_2$.

Let $n_3$ be the least natural number for which the inequality

$$\frac{[1 \cdot x_1] + \ldots + [n_1 x_1] + [(n_1 + 1)x_2] + \ldots + [n_2 x_2] + [(n_2 + 1)x_1] + \ldots + [n_3 x_1]}{\frac{n_3(n_3 + 1)}{2}} < x_1 + \varepsilon$$

is true. By Lemma 2, this number exists. Then we set $\alpha_{n_2+1} = \alpha_{n_2+2} = \ldots = \alpha_{n_3} = x_1$, etc.

Thus, we construct a sequence $(w_n)$, where

$$w_n = \frac{[1 \cdot \alpha_1] + [2 \cdot \alpha_2] + \ldots + [n \cdot \alpha_n]}{\frac{n(n + 1)}{2}}.$$

**Lemma 3.** *The sequence $(w_n)$ is convergent, i.e.,*

$$\lim_{n \to \infty} \frac{[1 \cdot \alpha_1] + [2 \cdot \alpha_2] + \ldots + [n \cdot \alpha_n]}{\frac{n(n + 1)}{2}}$$

*does not exist.*

**Proof.** Assume that $\lim_{n \to \infty} w_n$ exists. Then there exists $\delta$ such that $0 < \delta < x_2 - x_1 - 2\varepsilon$. By the Cauchy criterion of convergence of a sequence, there exists $n \in N$ such that, for any $m, l > n$, the inequality

$$|w_m - w_l| < \delta$$

is true.

By the construction of the sequence $(w_n)$, we have $w_{n_k} > x_2 - \varepsilon$ for odd $k$ and $w_{n_k} < x_1 + \varepsilon$ for even $k$. As usual, there exists $k \in N$ such that $n_k > n$. Then

$$|w_{n_{k+1}} - w_{n_k}| > x_2 - \varepsilon - (x_1 + \varepsilon) = x_2 - x_1 - 2\varepsilon > \delta,$$

i.e., $|w_{n_{k+1}} - w_{n_k}| > \delta$. We arrive at a contradiction. Thus, the limit $\lim_{n \to \infty} w_n$ does not exist.

We choose $x_1$ and $x_2$ such that $0 < x_1 < x_2$ and

$$\theta - 1 < x_i < \frac{2-\theta}{2}, \quad i = 1, 2,$$

and set $y_i = x_i - 1 + \theta > 0$ and $z_i = 2 - 2x_i - \theta$. Then $x_i + y_i + z_i = 1$ and

$$y_i + 2z_i = 2 - 2x_i - \theta + 2(2x_i - 1 + \theta) = \theta.$$

Let $\alpha_i = \tilde{\alpha}_i\{x_1, x_2\}$, $\beta_i = 2 - 2\alpha_i - \theta$, and $\gamma_i = \alpha_i - 1 + \theta$. We denote the lengths of the sequences formed by the digits "0," "1," and "2" in the number $x^*$ by $a_{i1} = [i \cdot \alpha_i]$, $a_{i2} = [i \cdot \beta_i]$, and $a_{i3} = [i \cdot \gamma_i]$, respectively. It is clear that

$$\alpha_i + \beta_i + \gamma_i = 1,$$

$$\beta_i + 2\gamma_i = \lambda.$$

We determine the mean value of the digits in the number $x^*$ as follows:

$$\frac{r_k}{k} = \frac{0 \cdot (a_{11} + \ldots + a_{(k-1)1}) + 1 \cdot (a_{12} + \ldots + a_{(k-1)2}) + 2 \cdot (a_{13} + \ldots + a_{(k-1)3}) + u_k}{a_{11} + \ldots + a_{(k-1)1} + a_{12} + \ldots + a_{(k-1)2} + a_{13} + \ldots + a_{(k-1)3} + v_k}$$

$$= \frac{[1 \cdot \beta_1] + \ldots + [(k-1) \cdot \beta_{k-1}] + 2[1 \cdot \gamma_1] + \ldots + 2[(k-1) \cdot \gamma_{k-1}] + u_k}{[1 \cdot \alpha_1] + [1 \cdot \beta_1] + [1 \cdot \gamma_1] + \ldots + [(k-1) \cdot \alpha_{k-1}] + [(k-1) \cdot \beta_{k-1}] + [(k-1) \cdot \gamma_{k-1}] + v_k},$$

where $u_k = n_1 + 2n_2$, $n_1$ is the number of digits "1," $n_2$ is the number of digits "2" in the $k$th block, and $v_k$ is the length of the $k$th block, i.e., the total number of digits "1," "1," and "2" in the $k$th block.

We estimate $u_k$ and $v_k$ as follows:

$$0 \leqslant u_k \leqslant [k \cdot \beta_k] + 2[k \cdot \gamma_k] \leqslant k \cdot \beta_k + 2k \cdot \gamma_k = k(\beta_k + 2\gamma_k) = k \cdot \theta,$$

$$0 \leqslant v_k \leqslant [k \cdot \alpha_k] + [k \cdot \beta_k] + [k \cdot \gamma_k] \leqslant k \cdot \alpha_k + k \cdot \beta_k + k \cdot \gamma_k = k(\alpha_k + \beta_k + \gamma_k) = k$$

and introduce the notation

$$A_k \equiv [1 \cdot \beta_1] + \ldots + [(k-1) \cdot \beta_{k-1}] + 2[1 \cdot \gamma_1] + \ldots + 2[(k-1) \cdot \gamma_{k-1}],$$

$$B_k \equiv [1 \cdot \alpha_1] + [1 \cdot \beta_1] + [1 \cdot \gamma_1] + \ldots + [(k-1) \cdot \alpha_{k-1}] + [(k-1) \cdot \beta_{k-1}] + [(k-1) \cdot \gamma_{k-1}].$$

According to properties of the integral part of a number, we find

(i) $\quad A_k \leqslant 1 \cdot \beta_1 + \ldots + (k-1) \cdot \beta_{k-1} + 2 \cdot 1 \cdot \gamma_1 + \ldots + 2 \cdot (k-1) \cdot \gamma_{k-1}$

$$= 1 \cdot (\beta_1 + 2\gamma_1) + \ldots + (k-1) \cdot (\beta_{k-1} + 2\gamma_{k-1}) = \theta(1 + 2 + \ldots + k - 1) = \theta \cdot \frac{(k-1)k}{2},$$

$$A_k > 1 \cdot \beta_1 - 1 + \ldots + (k-1) \cdot \beta_{k-1} - 1 + 2 \cdot 1 \cdot \gamma_1 - 2 + \ldots + 2 \cdot (k-1) \cdot \gamma_{k-1} - 2$$

$$= 1 \cdot (\beta_1 + 2\gamma_1) - 3 + \ldots + (k-1) \cdot (\beta_{k-1} + 2\gamma_{k-1}) - 3$$

$$= \theta(1 + 2 + \ldots + k - 1) - 3(k-1) = \theta \cdot \frac{(k-1)k}{2} - 3(k-1);$$

(ii) $$B_k \leqslant 1 \cdot \alpha_1 + 1 \cdot \beta_1 + 1 \cdot \gamma_1 + \ldots + (k-1) \cdot \alpha_{k-1} + (k-1) \cdot \beta_{k-1} + (k-1) \cdot \gamma_{k-1}$$

$$= 1 \cdot (\alpha_1 + \beta_1 + \gamma_1) + \ldots + (k-1) \cdot (\alpha_{k-1} + \beta_{k-1} + \gamma_{k-1}) = 1 + \ldots + k - 1 = \frac{(k-1)k}{2},$$

$$B_k > 1 \cdot \alpha_1 - 1 + 1 \cdot \beta_1 - 1 + 1 \cdot \gamma_1$$

$$- 1 + \ldots + (k-1)\alpha_{k-1} - 1 + (k-1)\beta_{k-1} - 1 + (k-1)\gamma_{k-1} - 1$$

$$= (\alpha_1 + \beta_1 + \gamma_1) - 3 + \ldots + (k-1)(\alpha_{k-1} + \beta_{k-1} + \gamma_{k-1}) - 3$$

$$= 1 + \ldots + k - 1 - 3(k-1) = \frac{(k-1)k}{2} - 3(k-1).$$

It follows from the last inequalities that

$$\frac{A_k}{\frac{(k-1)k}{2}} \to \theta \quad \text{and} \quad \frac{B_k}{\frac{(k-1)k}{2}} \to 1 \quad \text{as} \quad k \to \infty.$$

Then the mean value of digits for the number $x^*$ is equal to

$$\frac{r_k}{k} = \frac{A_k + u_k}{B_k + v_k} = \frac{\frac{A}{(k-1)} + \frac{u_k}{(k-1)k}}{\frac{B_k}{(k-1)} + \frac{v_k}{(k-1)k}} \to \frac{\theta + 0}{1 + 0} = \theta \quad \text{as} \quad k \to \infty.$$

We now show that the frequency of zero $\nu_0(x^*)$ does not exist. Indeed, if the indicated frequency exists, then it should be equal to

$$\lim_{k \to \infty} \frac{[1 \cdot \alpha_1] + \ldots + [(k-1) \cdot \alpha_{k-1}]}{B_k} = \lim_{k \to \infty} \frac{[1 \cdot \alpha_1] + \ldots + [(k-1) \cdot \alpha_{k-1}]}{\frac{(k-1)k}{2}}.$$

However, by Lemma 3, the last limit does not exist. Therefore, the frequency $\nu_0(x^*)$ also does not exist.

**Theorem 3.** *The set of numbers $W$ for which the frequencies of digits do not exist and the asymptotic mean value of digits is equal to a given number from $(0; 2)$ is a continual everywhere dense subset of the segment $[0; 1]$.*

***Proof.*** *Continuality.* We now show that different pairs $x_1$, $x_2$ and $x'_1$, $x'_2$ $(x_i \neq x'_j \; \forall i, j \in \{1, 2\})$ correspond to different numbers $x^*$ and $x'^*$. Assume the contrary, i.e., let $x^* = x'^*$. Then the $k$th series of zeros for the numbers $x^*$ and $x'^*$ are equal, i.e.,

$$a_{k1} = [k \cdot \alpha_k] = a'_{k1} = [k \cdot \alpha'_k]$$

for any $k \in N$. Thus,

$$|\alpha_k \cdot k - \alpha'_k \cdot k| < 1$$

for all $k \in N$, and, hence,

$$k \min_{i,j \in \{1,2\}} |x_i - x'_j| < 1$$

for all $k \in N$. However, the last inequality is not true for sufficiently large $k$. We arrive at a contradiction.

Since the set of pairs $(x_1, x_2)$ of numbers $x_1, x_2 \in \left(\theta - 1; \dfrac{2-\theta}{2}\right)$, $0 < x_1 < x_2$, is continual, the set of numbers for which the frequencies of digits do not exist and the asymptotic mean value of digits is equal to a given number from $(0; 2)$ is continual.

We now prove that the set $W$ is *everywhere dense*. Let

$$\left[\Delta^3_{\alpha_1 \alpha_2 \ldots \alpha_k (0)}; \Delta^3_{\alpha_1 \alpha_2 \ldots \alpha_k (2)}\right] \subset [0; 1],$$

be a certain cylindrical segment and let $x^* = \Delta^3_{\beta_1 \beta_2 \ldots \beta_n \ldots}$ be a number with asymptotic mean value of the digits but without any frequency of ternary digits, $\alpha_j, \beta_j \in \{0, 1, 2\}$. Since the asymptotic mean of the digits is independent of any finite number of the first digits in the ternary representation of a number, the number $x'^* = \Delta^3_{\alpha_1 \alpha_2 \ldots \alpha_k \beta_1 \beta_2 \ldots \beta_n \ldots}$ also has the asymptotic mean of the digits but does not have any frequency of the ternary digits. However, $x'^* \in \left[\Delta^3_{\alpha_1 \alpha_2 \ldots \alpha_k (0)}; \Delta^3_{\alpha_1 \alpha_2 \ldots \alpha_k (2)}\right]$. This implies that the set of numbers for which the frequencies of the digits do not exist but the asymptotic mean value of the digits is equal to a given number from $(0; 2)$ is an everywhere dense set in the segment $[0; 1]$.

**Theorem 4.** *The function $r$ of the asymptotic mean of the digits of a number has the following properties:*

1. *It is defined on an everywhere dense set in the segment $[0; 1]$.*

2. *It is not defined in an everywhere dense subset of the segment $[0; 1]$.*

3. *It takes all values from the set $[0; s - 1]$.*

4. *It has only one level of the full Lebesgue measure.*

***Proof.*** The first two properties follow from the facts established earlier.

Indeed, since the asymptotic mean value of the digits is independent of any finite number of the first digits in the ternary representation of a number, the set of points in which the function of frequency is defined is everywhere dense in $[0; 1]$.

There are numbers without asymptotic mean value of the digits. Thus, this is, e.g., the number

$$x^* = \Delta^s_{01001100001111\ldots \underbrace{0\ldots 0}_{2^n} \underbrace{1\ldots 1}_{2^n} \ldots}.$$

Since the asymptotic mean value of the digits is independent of any finite number of the first digits of a ternary number, the set of points in which the function of frequency is not defined is also an everywhere dense subset of $[0; 1]$.

3. Since

$$0 \leqslant \frac{1}{n} \sum_{i=1}^{n} \alpha_i(x) \leqslant s - 1,$$

the function $r(x)$ cannot take values higher than $s - 1$ and lower than 0. According to the algorithm presented above for the construction of a number with prescribed frequencies of the digits, the function $r(x)$ may take all values from $[0; s - 1]$.

4. By Theorem 1, only normal numbers have the property $r(x) = \frac{s-1}{2}$. Hence, the set of points at which the function $r$ takes the same fixed value $\frac{s-1}{2}$ is a set of full Lebesgue measure, i.e., the function $r$ has only one level of the full Lebesgue measure.

## REFERENCES

1. S. Albeverio, M. Pratsiovytyi, and G. Torbin, "Singular probability distributions and fractal properties of sets of real numbers defined by the asymptotic frequencies of their $s$-adic digits," *Ukr. Math. J.*, **57**, No 9, 1361–1370 (2005).
2. S. Albeverio, M. Pratsiovytyi, and G. Torbin, "Topological and fractal properties of real numbers which are not normal," *Bull. Sci. Math.*, **129**, No. 8, 615–630 (2005).
3. A. S. Besicovitch, "Sets of fractional dimension. 2. On the sum of digits of real numbers represented in the dyadic system," *Math. Ann.*, **110**, No 3, 321–330 (1934).
4. É. Borel, "Les probabilites denombrables et leurs applications arithmetiques," *Rend. Circ. Mat. Palermo*, **27**, 247–271 (1909).
5. H. G. Eggleston, "The fractional dimension of a set defined by decimal properties," *Quart. J. Math.*, Oxford Ser. 20, 31–36 (1949).
6. L. Olsen, "Normal and non-normal points of self-similar sets and divergence points of self-similar measures," *J. London Math. Soc.*, **2(67)**, No. 1, 103–122 (2003).
7. M. V. Prats'ovytyi, *Fractal Approach to the Investigation of Singular Distributions* [in Ukrainian], National Pedagogic University, Kyiv (1998).
8. M. V. Prats'ovytyi and H. M. Torbin, "Superfractality of the set of numbers having no frequency of $n$-adic digits, and fractal probability distributions," *Ukr. Mat. Zh.*, **47**, No. 7, 971–975 (1995); **English translation:** *Ukr. Math. J.*, **47**, No. 7, 1113–1118 (1995).
9. H. M. Torbin, "Frequency characteristics of normal numbers in different number systems," *Fract. Anal. Sumizh. Pyt.*, No 1, 53–55 (1998).
10. A. F. Turbin and N. V. Prats'ovytyi, *Fractal Sets, Functions, and Distributions* [in Russian], Naukova Dumka, Kiev (1992).